\documentclass[12pt]{article}
\usepackage[utf8]{inputenc}
\usepackage{amsmath}
\usepackage{amsfonts}
\usepackage{amssymb}
\usepackage{epsfig}
\usepackage{epstopdf}
\usepackage{indentfirst}
\usepackage{geometry}
\usepackage{graphicx}
\usepackage{caption}
\usepackage{subfig}
\usepackage[countmax]{subfloat}
\usepackage{rotating}
\usepackage{hyperref}
\usepackage{lscape}
\usepackage[numbers,sort&compress]{natbib}
\usepackage{tabularx}
\usepackage{setspace}
\usepackage{enumerate}
\usepackage{lineno}
\usepackage[section]{placeins}

\newtheorem{theorem}{Theorem}[section]
\newtheorem{lemma}{Lemma}[section]

\numberwithin{equation}{section}
\newcommand{\bet}{{\boldsymbol \eta}}

\newcommand{\blambda}{{\boldsymbol \lambda}}
\newcommand{\bepsilon}{{\boldsymbol \epsilon}}

\newcommand{\bzet}{{\boldsymbol \zeta}}
\def\bbR{{\mathbb R}}

\def\by{{\boldsymbol y}}

\def\bA{{\boldsymbol A}}
\def\bB{{\boldsymbol B}}
\def\bC{{\boldsymbol C}}

\def\ba{{\boldsymbol a}}
\def\bx{{\boldsymbol x}}

\def\bC{{\boldsymbol C}}
\def\bX{{\boldsymbol X}}

\def\bz{{\boldsymbol z}}

\def\bv{{\boldsymbol v}}

\def\btau{{\boldsymbol \tau}}

\def\bxi{{\boldsymbol \xi}}

\def\bbI{{\mathbb I}}

\def\cF{{\mathcal F}}

\def\cK{{\mathcal K}}
\def\cB{{\mathcal B}}
\def\cL{{\mathcal L}}

\def\cM{{\mathcal M}}

\graphicspath{%
    {converted_graphics/}% inserted by PCTeX
    {/}% inserted by PCTeX
}

\begin{document}

\title{Ill-posed linear inverse problems with box\\ 
constraints:  A new convex optimization approach}
\author{Henryk Gzyl \\
Center for Finance, IESA, Caracas.\\
{\small henryk.gzyl@iesa.edu.ve}
} 

\date{}
 \maketitle

%linenumbers  
\baselineskip=1.5 \baselineskip \setlength{\textwidth}{6in}
%\newpage

\begin{abstract}
Consider the linear equation $\bA\bx=\by,$ where $\bA$ is a $k\times N-$matrix, $\bx\in\cK\subset \bbR^N$ and $\by\in\bbR^M$ a given vector. When $\cK$ is a convex set and $M\not= N$ this is a typical ill-posed, linear inverse problem with convex constraints. Here we propose a new way to solve this problem when $\cK=\prod_j[a_j,b_j]$. It consists of regarding $\bA\bx=\by$ as the constraint of a convex minimization problem, in which the objective (cost) function is the dual of a moment generating function. This leads to a nice minimization problem and some interesting comparison results. More importantly, the method provides a solution that lies in the interior of the constraint set $\cK.$ We also analyze the dependence of the solution on the data and relate it to the Le Chatellier principle.
\end{abstract}


\begin{spacing}{0.1}
   \tableofcontents
\end{spacing}

\section{Introduction and preliminaries}
To establish the necessary notations we begin by stating the problem to be solved. Let $\bA$ be a 
$M\times N-$matrix and consider the linear equation:
\begin{equation}\label{prob1}
 \begin{aligned}
&\bA\bx=\by\\
\mbox{With}\;\;\; \bx\in\prod_{j=1}^N [a_j,b_j]& \subset\bbR^N,\;\;\;\mbox{and}\;\;\by\in\bbR^M.
\end{aligned}
\end{equation}
In many problems $M\not= N$ and when equal, the matrix might not be invertible. This is a typical ill-posed, linear inverse problem, and because of the constraint $\bx\in\prod_{j=1}^N [a_j,b_j]$ placed upon the solution, it is said to have convex constraints. 

Problem \ref{prob1} admits the important extension that is essential in many applications. Consider \cite{Go} as an example of an over-determined problem or \cite{Var1} as an example of an under-determined problem. Ill-posed problems have been studied since long ago and many techniques to solve them have been proposed. Consider for example \cite{VL} for the use of maximum likelihood techniques or \cite{BB2} and \cite{BB}, for an update on regularization techniques and \cite{Gz} for an application of maximum entropy in the mean to determine joint probabilities from marginals. The list of problems in applied mathematics, that leads to problems like \eqref{prob1} or \eqref{prob2}, is very long.  

The method that we propose is a notational extension of the problem that occurs in actual practice, which can be stated as: Determine $\bx$ and $\bepsilon$ that satisfy:

\begin{equation}\label{prob2}
 \begin{aligned}
&\bA\bx + \bepsilon=\by\\
\mbox{Such that}\;\;\; \bx\in\prod_{j=1}^N [a_j,b_j] \subset\bbR^N,\;\;&\bepsilon \in\prod_{j=1}^K [c_j,d_j] \subset\bbR^M\;\;\;\mbox{and}\;\;\by\in\bbR^M.
\end{aligned}
\end{equation}
This describes the case in which the response $\by$ is measured up to an additive noise $\bepsilon$ that can be modeled as a quantity to be determined in $\prod_{j=1}^K [c_j,d_j].$ The values $c_j,d_j$ are to be determined by the modeler from the statistical information about the measurement procedure.

The method to solve Problems \ref{prob1} and \ref{prob2} that we propose, is reminiscent of the method of maximum entropy in the mean. Its consists of transforming \eqref{prob1} it into a variational problem, consisting of the minimization of a Bregman divergence subject to appropriate constraints. The similarity to the maxentropic procedure lies in the way of defining the divergence as the Fenchel-Lagrange dual of an appropriate Laplace transform. The divergence is defined on the constraint set and it avoids the intermediate step of obtaining the solution as an expected value with respect to an unknown probability that is obtained by an entropy maximization method.

Having solved Problem \ref{prob1}, the extension to solve problem \ref{prob2} involves a slight change in notations. To address \ref{prob1}, consider the measure space $(\Omega,\cF,Q)$ where $\Omega=\prod_{j=1}^N [a_j,b_j]$ is the set of constraints,  $\cF=\cB(\Omega)$ denotes the Borel subsets of $\Omega.$ The (reference) measure $Q$ is any measure such that the closed convex hull of its support equals $\Omega.$ We shall consider

\begin{equation}\label{ref}
dQ(\bxi) = \otimes_{j=1}^N\big(\epsilon_{a_j}(d\xi_j) + \epsilon_{a_j}(d\xi_j)\big).
\end{equation}
The notation $\epsilon_{a}(d\xi)$ stands for the point mass (Dirac) measure at the point $a.$ This reference measure is chosen for several reasons: First, its joint Laplace transform is trivial to compute, second, because its convex support is $\Omega$ and third, because of it leads to easy and interesting computations. 

Let $\bX:\Omega\to\Omega$ defined by $\bX(\bxi)=\bxi$ denote the identity (or coordinate) mapping. The (multi-dimensional Laplace transform) of $Q$ is given by:

\begin{equation}\label{LT}
\bzet({\btau}) = E_Q[e^{\langle\btau,\bX\rangle}] =  \prod_{j=1}^N\zeta_j(\tau_j) = \prod_{j=1}^N\big(e^{a_j\tau_j}+e^{b_j\tau_j}\big),\;\;\;\btau\in\bbR^N.
\end{equation}
In the probability/statistics literature the function $M(\tau)$ defined right below is called the moment generating function.
\begin{equation}\label{momgf}
M(\btau)=\ln\bzet(\btau) = \sum_{j=1}^NM_j(\tau_j) = \sum_{j=1}^N\ln\big(e^{a_j\tau_j}+e^{b_j\tau_j}\big).
\end{equation}
The following lemma is used to obtain the Legendre-Fenchel dual of $M(\btau).$ Since the variables $\tau_j$ appear separated in $M(\btau)$ we shall drop reference to the $j-$th subscript and put it back in when needed.

\begin{lemma}\label{lem1}
Let $\bxi\in (a,b).$ The equation $\bxi=\nabla_\tau M(\tau)$ has a unique solution that establishes a bijection between $\bbR$ and $(a,b).$
\end{lemma}
The proof is computational. Note that
$$\xi = \frac{\partial M(\btau)}{\partial \tau} = a\frac{e^{a\tau}}{e^{a\tau}+e^{b\tau}}+b\frac{e^{b\tau}}{e^{a\tau}+e^{b\tau}}$$
has a solution given by
\begin{equation}\label{chvar}
e^{\tau} = \bigg(\frac{\xi-a}{b-\xi}\bigg)^{1/D}\;\;\;\Leftrightarrow\;\;\;\tau=\frac{1}{D}\ln\bigg(\frac{\xi-a}{b-\xi}\bigg).\;\;\; D=b-a.
\end{equation}
Note as well that $\xi(\tau) \to b$ (respectively to $a$) whenever $\tau\to-\infty$) (respectively to $+\infty$).
\begin{lemma}\label{lem2}
The Legendre-Fenchel dual of $M(\tau)=\ln\big(e^{a\tau}+e^{b\tau}\big)$ is defined to be $\psi(\xi)=\sup\{\xi\tau-m(\tau): \tau\in\bbR\}.$ It is
\begin{equation}\label{LF1}
\psi(\xi) = \frac{\xi-a}{D}\ln\big(\frac{\xi-a}{D}\big) + \frac{b-\xi}{D}\ln\big(\frac{b-\xi}{D}\big).
\end{equation}
To finish, the Legendre-Fenchel dual of $M(\tau)$ is given by

\begin{equation}\label{LF2}
\Psi(\bxi) = \sum_{j=1}^N\frac{\xi_j-a_j}{D_j}\ln\big(\frac{\xi_j-a_j}{D_j}\big) + \frac{b_j-\xi_j}{D_j}\ln\big(\frac{b_j-\xi_j}{D_j}\big).
\end{equation}
We also have:
\begin{equation}\label{deriv}
\big(\Psi'\big)^{-1}(\xi) = M'(\tau)\;\;\;\mbox{or equivalently}\;\;\; \Psi'(\xi) = \big(M'\big)^{-1}(\tau).
\end{equation}
\end{lemma}
Clearly, each $\psi(\xi_j)$ and $\Psi(\bxi)$ are convex on $\Omega.$ \\
The proof uses the results in Lemma \ref{lem1} and simple arithmetic to obtain (\ref{LF1}), and putting the subscripts back in, we obtain (\ref{LF2}). As usually it is consistent to take $0\ln 0=0.$ We mention that (\ref{deriv}) is a standard identity relating the derivatives of a function and that of its Legendre-Fenchel dual.

\subsection{The contents of this paper} 

Notice that the Hessian of $\Psi(\bxi)$ is a sum of strictly convex functions in each variable. In Section 2 we consider some properties of the Bregman divergence associated with $\Psi(\bx).$ In particular, we will obtain some interesting comparison results.

In Section 3 we shall render problem (\ref{prob1}) as a variational problem consisting of minimizing $\Psi(\bx)$ subject to the equation $\bA\bx=\by$ as a constraint. In Section 4 we consider a discrete version of the problem dealt with in \cite{C}, 

Section 4 is devoted to the analysis of the dependence of the solution on the data vector $\by.$ We prove that some special condition on $\bA,$ the Jacobian $\partial(\blambda)/\partial(\by), where$$\blambda$ is the Lagrange multiplier introduced in Section 3. This will determine the dependence of the optimal solution $\bx^*$ of \eqref{prob1} on the data. In section 5 we consider an application to a discrete vesion of a problem dealt with in \cite{C}, consisting of determining a bounded function $f(j): j=1,...,N$ (the initial data) if the final data $g(i): i=1,...,K$ is known at a number of points in the state space.

\section{A divergence derived from $\Psi$ and a comparison result}
Since $\Psi(\bxi)$ is a strictly convex function defined on $\Omega$ we can associate to it a Bregman divergence defined by: 

\begin{theorem}\label{breg0}
The Bregman divergence determined by $\Psi$ is defined by

$$\delta_\Psi^2(\bxi,\bet) = \Psi(\bxi)-\Psi(\bet)-\langle(\bxi-\bet,\nabla\Psi(\bet)\rangle.$$
For $\Psi(\bxi)$ as in (\ref{LF2}) we have
\begin{equation}\label{breg1}
\delta_\Psi^2(\bxi,\bet)= \sum_{j=1}^N\left[\frac{\xi_j-a_j}{D_j}\ln\bigg(\frac{\xi_j-a_j}{\eta_j-a_j}\bigg)+\frac{b_j-\xi_j}{D_j}\ln\bigg(\frac{b_j-\xi_j}{b_j-\eta_j}\bigg)\right].
\end{equation}
 Furthermore
\begin{equation}\label{propberg}
\delta_\Psi^2(\bxi,\bet) \geq 0 \;\;\;\mbox{and} \;\;\;\delta_\Psi^2(\bxi,\bet)=0 \Leftrightarrow \bxi=\bet.
\end{equation}
\end{theorem}
As $\Psi(\bxi)$ is reminiscent of an entropy function so is $\delta_\Psi^2(\bxi,\bet)$ reminiscent of a cross-entropy (or a Kullback-Leibler divergence) defined on $\prod[a_j,b_j]$ if we think of $p(\xi_j)=(\xi_j-a_j)/(b_j-a_j)$ as the probability of picking a point uniformly in $(a_j,\xi_j)\subset (a_j,b_j).$ 

\subsection{A comparison result}
Notice that since the Hessian of $\Psi$ is a diagonal matrix, we have:
$$\delta_\Psi^2(\bxi,\bet) = \Psi(\bxi)-\Psi(\bet)-\langle(\bxi-\bet),\nabla\Psi(\bet)\rangle$$
$$ = \sum_{j=1}^N\int_{\xi_j}^{\eta_j}(\eta_j-x_j)\frac{1}{(b_j-a_j)}\bigg(\frac{dx_j}{(x_j-a_j)}+\frac{dx_j}{(b_j-x_j)}\bigg).$$
Drop the subscript and notice that $(x-a)(b-x)$ reaches its maximum value $(b-a)^2/4$ at $(a+b)/2.$ Therefore
$$
\begin{aligned}
\delta_{\Psi}^2(\bxi,\bet) = \sum_{j=1}^N\int_{\xi_j}^{\eta_j}(\eta_j-x_j)\frac{dx_j}{(x_j-a_j)(b_j-x_j)}\\
\geq \sum_{j=1}^N\int_{\xi_j}^{\eta_j}\frac{4}{(b_j-a_j)^2}(\eta_j-x_j)dx_j = 2\sum_{j=1}^n\bigg(\frac{(\eta_j-\xi_j)}{(b_j-a_j)}\bigg)^2.
\end{aligned}
$$
These calculations prove the following result.
\begin{theorem}\label{comp1}
With the notations introduced above, for any $\bxi,\bet \in \Omega$ we have:
\begin{equation}\label{comp1.1}
\delta_{\Psi}^2(\bxi,\bet) = 2\sum_{j=1}^N\bigg(\frac{(\eta_j-\xi_j)}{(b_j-a_j)}\bigg)^2.
\end{equation}
\end{theorem}
To obtain an upper bound, again drop the subscripts and notice that
$$
\int_{\xi}^{\eta}(\eta-x)\frac{1}{(b-a)}\bigg(\frac{dx}{(x-a)}+\frac{dx}{(b-x)}\bigg) \leq \frac{(\eta-\xi)}{(b-a)}\int_\xi^\eta d\ln\bigg(\frac{a-x}{b-x}\bigg).
$$
Thus we have:
\begin{theorem}\label{comp2}
With the notations introduced above, for any $\bxi,\bet \in \Omega$ we have:
\begin{equation}\label{comp2.1}
\delta_{\Psi}^2(\bxi,\bet) \leq \sum_{j=1}^N\bigg(\frac{(\eta_j-\xi_j)}{(b_j-a_j)}\bigg)\big(\phi(\eta_j)-\phi(\xi_j)\big)
\end{equation}
where we put
$$\phi(x) = \ln\bigg(\frac{x-a}{b-x}\bigg).$$
\end{theorem}

\section{The solution to problems \ref{prob1}-\ref{prob2}}
The solution of the two problems is the same. We work out the solution to Problem \ref{prob1} and explain how to modify the notations to obtain the solution to \ref{prob2}.

\subsection{The solution to problem \ref{prob1}}
Here we make use of the dual moment generating function $\Psi(\bxi)$ and its associated Bregman divergence $\delta_\Psi^2(\bxi,\bet)$ to solve problem (\ref{prob1}). We already mentioned that $\Psi(\bxi)$ can be regarded as an entropy and $\delta_\Psi^2(\bxi,\bet)$  as a cross entropy. Instead of solving (\ref{prob2}), we propose the alternative version

\begin{equation}\label{prob3}
\mbox{Determine}\;\;\;\; \bxi^* \;\;\;\;\mbox{that realizes}\;\;\;\; \min\{\Psi(\bxi)| \bA\bxi=\by\}.\;\;\;\;\;\;\;\end{equation}

The standard argument invoking Lagrange multipliers and the result in (\ref{deriv}) implies that the solution to (\ref{prob3}) is given by:
\begin{equation}\label{repsol1}
\bxi_j^*= \frac{a_je^{a_j(\bA^*\blambda^*)_j} + b_je^{b_j(\bA^*\blambda^*)_j}}{e^{a_j(\bA^*\blambda^*)_j}+e^{a_j(\bA^*\blambda^*)_j}},\;\;\;\;\;\;j = 1,...N..
\end{equation}
Here $\bA^*$ denotes the transpose of $\bA$ and $\blambda^*$ is a Lagrange multiplier that has yet to be determined. To determine $\blambda^*$ we use the following variant of a duality argument.\\
Consider the surface in $\cM$  $\blambda\to\bet(\blambda)$  parameterized by $\blambda\in\bbR^M$ given by:
\begin{equation}\label{surf}
\bet(\blambda)_j= \frac{a_je^{a_j(\bA^*\blambda)_j} + b_je^{b_j(\bA^*\blambda)_j}}{e^{a_j(\bA\blambda)_j}+e^{a_j(\bA^*\blambda)_j}},\;\;\;j = 1...N.
\end{equation}
That is, $\bx^*$ in (\ref{repsol1}) might turn out to be a point in this surface if we find the $\blambda^*.$ It is a simple calculation to verify that if $\bxi$ is any possible solution to (\ref{prob3}), then for any $\blambda\in\bbR^M$ we have:
$$\delta_\Psi^2(\bxi,\bet(\blambda)) \geq 0 \Leftrightarrow \Psi(\bxi) \geq \langle\blambda,\by\rangle - M(\bA^*\blambda).$$
Recall that $M$ is the moment generating function introduced in \eqref{momgf}. To finish, note that the first order condition that determines the maximizer of the ``dual entropy'':
\begin{equation}\label{dualent}
\Sigma(\by,\blambda)=\langle\blambda,\by\rangle - M(\bA^*\blambda)
\end{equation}
asserts that the right hand side of (\ref{repsol1}) is indeed a point in the surface \eqref{surf}. The interest in this argument lies in the fact that in order to solve (\ref{prob3}) it suffices to maximize a strictly convex function $\Sigma(\by,\blambda)$ over $\bbR^M.$ We collect  all these comments and computation under:

\begin{theorem}\label{solinv}
With the notations introduced above, the solution to the ill-posed linear inverse problem with box constraints (\ref{prob3}) is given by (\ref{repsol1}) in which the Lagrange multiplier $\blambda^*$ is obtained as the maximizer of the dual entropy (\ref{dualent}). Furthermore
\begin{equation}\label{dual1}
\Psi(\bxi^*) = \langle\blambda^*,\by\rangle - M(\bA^*\blambda^*) = \Sigma(\by,\blambda^*).
\end{equation}
\end{theorem}

This is a standard result in duality theory. See \cite{L} for example. Theorem \ref{solinv} leads to a simple numerical procedure to obtain $\blambda^*.$  We add that the $\blambda^*$ obtained minimizing \eqref{dualent}, is the same as the Lagrange multiplier that is obtained minimizing $\Psi(\bx)$ subject to \eqref{prob1} as constraint using the Lagrange procedure. This fact is used below when invoking the Le Chatelier principle to analyze the dependence on the initial data.

\subsection{The solution to Problem \ref{prob2}}
The adaptation is rather straight forward. Let us write \eqref{prob2} as
\begin{equation}\label{prob3}
\begin{split}
\bB\bz& = \by\\
\mbox{where}\;\;\;\bB=[\bA\; \bbI]\;\;\;&\mbox{and}\;\;\;\bz={\bx\atopwithdelims( )\bepsilon}.
\end{split}
\end{equation}

Now $\bB$ is an $(N+M)\times M-$matrix, $\bbI$ is the $M\times M-$identity matrix, $\bz\in\bbR^{N+M},$ and the constraint set is $\prod_{j=1}^N[a_j,b_j]\times\prod_{j=1}^M[c_j,d_j].$ The reference measure is now

\begin{equation}\label{ref1}
dQ(\bxi) = \otimes_{j=1}^N\big(\epsilon_{a_j}(d\xi_j) + \epsilon_{a_j}(d\xi_j)\big)\otimes_{j=1}^M\big(\epsilon_{c_j}(d\eta_j) + \epsilon_{a_j}(d\eta_j)\big).
\end{equation}

To complete the description of the building blocks, the moment generating function is:

\begin{equation}\label{momgf1}
M(\btau) = \sum_{j=1}^N\ln\big(e^{a_j\tau_j}+e^{b_j\tau_j}\big) + \sum_{j=1}^M\ln\big(e^{c_j\tau_{j+N}}+e^{d_j\tau_{N+j}}\big).
\end{equation}

Repeating the procedure that leads to \eqref{repsol1} we obtain:

\begin{equation}\label{repsol2}
\bxi_j^*= \frac{a_je^{a_j(\bA^*\blambda^*)_j} + b_je^{b_j(\bA^*\blambda^*)_j}}{e^{a_j(\bA^*\blambda^*)_j}+e^{a_j(\bA^*\blambda^*)_j}}\;\;\;\mbox{and}\;
\;\;\bepsilon_j=\frac{c_je^{c_j\blambda^*_j} + d_je^{d_j\blambda^*_j}}{e^{c_j \blambda^*_j}+e^{d_j \blambda^*_j}}.
\end{equation}

This time, the optimal $\blambda^*$ has to be found minimizing the convex function:

\begin{equation}\label{dualent1}
\Sigma(\by,\blambda)=\langle\blambda,\by\rangle - M(\bB^*\blambda)
\end{equation}
where now $M$ is given by \eqref{momgf1} in which $B^*$ denotes the transpose of $\bB.$

\section{Dependence on the data and the Le Chatelier principle}
Clearly, the solution $\bx^*$ given by either \eqref{repsol1} or \eqref{repsol2} depends smoothly on $\blambda^*.$ So if we establish the dependence of $\blambda^*$ on the data vector $\by,$ we can determine the dependence of $\bx^*$ on $\by.$  

Notational convention: in order to unclutter the typography, we drop the $^*$ denoting optimality and keep the detail in mind.

A simple computation shows that the Condition for $\blambda$ to be a minimizer of \eqref{dualent} is that $\bA\bx(\blambda)=\by.$ Since
$$\frac{\partial x_j}{\partial y_k} = \sum_{m=1}^M \frac{\partial x_j}{\partial \lambda_m}\frac{\partial \lambda_m}{\partial y_k} .$$ 
 
So, the essential consists of computing $\frac{\partial x_j}{\partial \lambda_m}.$ The computation is simple and it leads to

\begin{equation}\label{interm1}
\frac{\partial x_j}{\partial \lambda_m}=-\sum_{m=1}^M A_{k,j}\bigg(\frac{(b_j - a_j)}{e^{(b_j-a_j)\tau_j/2} + e^{-(b_j-a_j)\tau_j/2}}\bigg)^2,
\end{equation}
where $\tau_j$ should be replaced by $(\bA^*\blambda)_j.$ Therefore, had we differentiated $\bA\bx=\by$ with respect to $y_k,$ and denoting the term inside the parentheses by $C_j,$ we would have obtained that the $i$th is given by:
\begin{equation}\label{interm2}
-\sum_{j=1}^N \sum_{m=1}^M A_{i,j}C^2_j A^*_{j,m}\frac{\partial \lambda_m}{\partial y_k} = \delta_{i,k}.
\end{equation}
Or in vector notation $-\bA\bC\bA^*\partial(\blambda)/\partial(\by)=\bbI$ where $\bC$ is the diagonal matrix described a few lines above, and  
$$\frac{\partial(\blambda)}{\partial(\by)}_{k,m} = \frac{\partial \lambda_m}{\partial y_k}$$
stands for the Jacobian matrix.
 
To sum up:

\begin{theorem}\label{sens}
With the notations introduced a few lines above, suppose that the matrix $\bA\bC\bA^*$ (which, to begin with, is positive semidefinite) is invertible, then \eqref{interm2} implies:
\begin{equation}\label{jac}
\frac{\partial(\blambda)}{\partial(\by)} = -\big(\bA\bC\bA^*\big)^{-1}.
\end{equation}
\end{theorem}
This asserts that $\blambda$ is continuous with respect to $\by,$ which taken jointly with \eqref{repsol1} implies that $\bx^*$ is continuous with respect to $\by.$ Now we can make use of Theorems \ref{comp1} or \ref{comp2} to relate the change of the the divergence between the solutions for different data with respect of the changes of the solutions.

Theorem is a weaker version of the Le Chatelier principle presented in \cite{S}. Here we follow the argument presented in \cite{EO} to prove an extended version of Theorem \ref{sens}.

 Let the Lagrangian function 
$\cL(\bx,\blambda,\by):\cK\times\bbR^N\times\bA(\cK)\to \bbR$ be defined by
$$\cL(\bx,\blambda,\by) = \Psi(\bx) + \blambda^t\big(\bA\bx-\by)$$
For $i=1,2,$ let $(\bx_i^*,\blambda_i^*)$ be the solutions to the problem on minimizing the Lagrangian when the data is $\by_i.$ Then:

\begin{gather}\nonumber
\cL(\bx^*_1,\blambda^*_1,\by_1)\;\;\leq\;\; \cL(\bx^*_2,\blambda^*_2,\by_1)\nonumber\\
\cL(\bx^*_2,\blambda,^*_2,\by_2)\;\; \leq\;\;\cL(\bx^*_1,\blambda^*_1,\by_2).\nonumber
\end{gather}
From this it readily drops out that:
$$\cL(\bx^*_2,\blambda^*_2,\by_2) - \cL(\bx^*_2,\blambda^*_2,\by_1) \leq \cL(\bx^*_2,\blambda^*_2,\by_1) - \cL(\bx^*_1,\blambda^*_1,\by_1).$$
Now, substituting the definition of each term, simplifying  and rearranging we obtain the following extension of \eqref{jac}:

 $$(\blambda_1^*-\blambda_2^*)^t(\by_1-\by_2) \leq 0.$$

As $\bx$ and $\by$ live in different spaces, one thing we can do to relate the change in $\bx^*$ to the change in $\by,$ is to replace $\by_i$ with $\bA\bx_i$ in the previous identity, and rewrite it as:
 
$$(\blambda_1^*-\blambda_2^*)^t(\bA\bx_1-\bA\bx_2) = (\bA^t\blambda_1^*-\bA^t\blambda_2^*)^t(\bx_1-\bx_2)\leq 0.$$

\section{Example}
As mentioned in the introduction, we consider a discrete version of a problem dealt with in 
\cite{C}, consisting of determining the initial values of some observable of a Markovian system (a diffusion, or a Brownian motion, to be specific) on $\bbR^d.$ If we denote the transition semigroup at some time $T>0$ by $P,$ the problem consists of determining the initial data $f$ from the knowledge of $g=Pf.$ It is solved in \cite{C} in a Hilbert space setting by a combination of spectral and regularization techniques.

The discrete version of that problem, cast as a variation on the theme of Problem
\eqref{prob1},  goes as follows:

Let $\bA=P$ be a Markov transition matrix of a process having discrete state space $\{1,2,...,N\}.$ Let $f:\{1,...,N\}$ be some unknown observable and suppose that after some time we can only observe its value, denote it by $g(i),$ at some subset of the state space. Denote it by $\{1,...,M\}$ to simplify the notations. The problem consists of determine a function $f$ such that:
\begin{equation}\label{prob2}
g(i) = \sum_{j=1}^N P_{i,j}f(j),\;\;\;i=1,...,M.
\end{equation}
The case of interest is when $N$ is very large and it is costly to observe the state of the system, so $M$ is considerably smaller than $N$ and the problem is quite undetermined. If the modeler requires the bounds on the initial datum to be $0\leq f(j)\leq b$ for all $j=1,...,N,$ then an application of \eqref{repsol1} yields:

\begin{equation}\label{repsola}
f^*(j) = b\frac{1}{1 + e^{b(P\lambda)_j}}.
\end{equation}
The dependence of $f^*$ on the data $g$ is as explained in the previous section.

\section{Final comments}
A small change in the notations can be used to deal with the case when the constraints set in \eqref{prob1} is a convex set like $\cK = \{\bx\in\bbR^N|\,\bx=\sum_{i=1}^K a_i\bv_i=d\},$ where the $\bv_i\in\bbR^N$ are the vertices of a polyhedron, and the vector $\ba\in\bbR^K$ satisfies some box constraints.

We mention that the method proposed above is a relative of the method of maximum entropy inn the mean. That method was established  in \cite{DG}. Our proposal differs from the maxentropic procedure in the fact that maxentropic solution is obtained as the mean value of appropriate random variables, that range over the set of constraints, with respect to a probability obtained by an application of the standard method of maximum entropy. In the procedure proposed above, there is no need of the determination of a probability because the objective function is defined directly on the constraint space. See the application of the method of maximum entropy in the mean presented in \cite{Gz} for comparison.

\section{Declaration: Conflicts of interest/Competing Interests}
I certify that there is no actual or potential conflict of interest in relation to this article.

 \end{document}